\documentclass[a4paper,10pt]{article}

\usepackage{amsfonts,amssymb,mathrsfs,amscd}
\def \R {{\mathbb {R}}}

\def \Z {{\mathbb {Z}}}

\def\uu{\bigsqcup}
\def\eps{\varepsilon}

\title{\bf Ergodic Theorems, Almost Invariant Sets, \\
           Value Distributions  of Time Averagings}
\author{\bf V.V. Ryzhikov}
\date{Lomonosov Moscow State University}
\begin{document}
\Large
\maketitle

 {\bf Abstract.} {\large Based on classical facts of ergodic theory, we demonstrate how deviations of ergodic averages from the constants to which they converge can be structured. The deviations tend to zero almost everywhere and with respect to the integral norm. They are asymptotically nearly invariant under the action due to averaging. In this situation, the question of the distribution of the values of such deviations is meaningful. We show that  these distributions can be weakly close to any given distribution if we apply a homothety to the latter  changing the scale on the value line. We  provide in laconic  manner  related proofs of  von Neumann, Birkhoff, Wiener theorems, and Rokhlin's lemma for $\Z^d$-actions, ensuring the text is self-contained.}

\section{Introduction}
Numerous ergodic theorems represent an independent, extensive field in operator theory and the theory of dynamical systems; see, for example, the books
\cite{K}, \cite{T}. Our aim is to answer the  question of what deviations in Birkhoff's theorem  might look like. To give a close presentation we do the following. Show that ergodic theorems can be used to prove Rokhlin's lemmas and show the application of such lemmas to proving the convergences  almost everywhere. The convergences almost everywhere, very briefly, is proved according to the following scheme. Averaging is performed over almost-invariant (Follner) subsets of the group. We assume that the function being averaged has zero mean and observe the convergence of the time averages to zero. If there is no convergence, then, due to the ergodicity of the action, we have a divergence uniform  of the deviations. We consider an extremely large
almost-invariant time set (a large tower), little different from the entire space, but endowed with the structure of a tower partition. In this large tower, due to the divergence of the means, we find an almost-invariant set on which
our averaged function has an integral bounded away from zero. Proceeding to the consideration of sequences of such invariant sets, by virtue of the weak ergodic theorem for our action, we have the convergence of such integrals to zero, thereby arriving at a contradiction with the divergence almost everywhere.

Then almost-everywhere convergence is used to demonstrate the following observation. Deviations from the spatial mean tend to zero almost everywhere. They are asymptotically almost invariant under the action due to averaging. In this situation, the question of the distribution of the values of such deviations is meaningful. It turns out that for any ergodic free action of the group $\Z^d$, these distributions can be weakly close to any predetermined distribution when the scale on the value line changes.

\section{ Convergence  in norms}
For an ergodic automorphism $T$ of the probability space $(X,\mu)$, a function $f\in L_1(X,\mu)$, the von Neumann and Birkhoff theorems assert the convergence in $L_1$ and,
respectivly,  almost everywhere to  $ \int_X f\, d\mu$ as $N\to\infty$
of time ergodic averages
$$P_N f(x):=\frac 1 N \sum_{i=1}^{N} f(T^ix). $$
Let us recall  ideas of the proofs of these theorems.

\vspace{3mm}
\bf Theorem 2.1. \it $\left\|P_N f-\int f d\mu\right \|_{1,2}\to 0.$ \rm

\vspace{3mm}
Proof.
Let $f\in L_\infty(\mu)$ have zero mean. We show that $\|P_Nf\|_2\to 0$.
Since $\|TP_N-P_N\|\leq 2/N\to 0$, we have
$$ \|TP^\ast_N P_N-P^\ast_NP_N\|\to 0.$$
If $P$ is a limit point for the sequence $P_N$, then
$$ TP^\ast P_N-P^\ast P=0, \ \ TP^\ast P_Nf =P^\ast Pf.$$
Since the automorphism $T$ is ergodic, the function $P^\ast Pf$ is constant, hence,
it is zero. From this we obtain $(P^\ast_NP_Nf\,,\,f)\to \|f\|^2_2\to 0.$
Since the functions $P_Nf$ converging to 0 in measure are bounded, we have $\|P_Nf\|_1\to 0$. The norms of the operators $P_N:L_p\to\L_p$ are equal to 1, therefore
$\|P_Nf\|\to 0, \ \ f\in L_1.$
The theorem is proved.

The same method  works gor group actions.
Denote
$$Q_N=\{z=(z_1,\dots,z_d)\, :\, 1\leq z_1,\dots,z_d\leq N\}.$$

\vspace{3mm}
\bf Theorem 2.2. \it Let $f\in L_1(\mu)$.  For an ergodic   $\{T^z\}$-action we have 
$\|P_Nf -\int_X f d\mu\|_{1,2} \to 0$,  where $P_N=N^{-d}\sum_{z\in Q_N}T^z$.  

\vspace{3mm}
The assertion about the equivalence of the convergences $P_N\to_s 0$ and
$TP^\ast_N P_N\to_w 0$ can be generalized in various ways;
we give an example.

\vspace{3mm}
\bf Theorem 2.3. \it Let $\{T_g\, :\, g\in G\}$ be an ergodic action by automorphisms of the probability space of a countable infinite group, and let $P_j=\sum_{g} w_j^g T_g$, where  $\sum_g w^g_j=1$ and  $ w_j^g\leq 0.$ The strong convergence  $P_j\to_s \Theta$,  $
\Theta f\equiv \int_X f d\mu$,   is equivalent to the weak convergence
$T_g(P^\ast_j P_j)^{2^k}\to_w \Theta, \ j\to\infty. $
\rm

\vspace{3mm}
The proof is not significantly different from the previous one.

\vspace{3mm}
\bf Discrete averaging over spherical fibers. \rm
We give an example of an ergodic theorem, leaving the proof to the reader

\vspace{3mm}
\bf Theorem 2.4. \it Let $F_j=\{z\in \Z^3\,:\, j<|z|<j+c\}$. If the constant $c>0$ is sufficiently large, then for an ergodic  $\Z^3$-action the  averages
$P_jf=|F_j|^{-1} \sum_{z\in F_j}T^zf$ 
converge  to  $\int_Xfd\mu$ in $L_1(\mu)$. \rm

\vspace{3mm}
For the proof, we must use the weight distribution of the operators $P^\ast_jP_j$.
For $c>2$, the distribution of these weights is relatively smooth, and no difficulties with average convergence arise. However, for $c<1$, the problem no longer appears simple.

Continuous $\R^n$-analogs of the theorem are considered in \cite{BR}. The author thanks Benjamin Weiss communicated to author on the paper \cite{FKW} containing an example of averaging over a time manifold of lower dimension, together with effective applications to Ramsey-type problems on measurable sets.

\section{Rokhlin lemma and Birkhoff theorem }
\bf Theorem 3.1. \it Let $T$ be an ergodic automorphism of the probability space $(X,\mu)$. For every $\eps>0$ and a natural number $n$, there exists a measurable set $B$ such that
$$X=\uu_{i=1}^{n} T^{i-1}B\uu E, \ \ \mu(E)<\eps.$$ \rm

\vspace{2mm}
The complement of $E$ is called the Rokhlin tower (or Rokhlin-Halmos tower, see \cite{KSF}) of height $n$.

Proof. For an arbitrary set $D$ of small positive measure, consider the Kakutani partition over $D$:
$$X=\uu_{h=1}^\infty\uu_{i=1}^{h} T^{i-1}D_h,$$
where $$D_h =\{x\in D: T^hx\in D, \ T^ix\notin D, \ 0<i<h\}.$$
Such sets $\uu_{i=1}^{h} T^{i-1}D_h$ are called columns, and $D_h$ is the base of the column.
Since the measure of the entire space is 1, we find $k> 1/\mu(D)$ for which $\mu(D_k)>0$ (the number $k$ can be chosen arbitrarily large, decreasing the measure of the set $D$).
Consider the Kakutani partition over $D_k$. Now all columns have a height no less than $k$.
The base $B$ of the desired tower $\uu_{i=1}^{n} T^{i-1}B$ is obtained as the union of all the floors of columns with numbers multiples of $n$. The complement of this tower has measure no greater than $n/h<\eps$. The lemma is proved.

In proving Birkhoff's theorem, we will use this lemma as follows:
a (very large) tower of height $n$ is the union of pieces of trajectories starting at the base $B$. Working with each piece separately, we obtain a globally measurable object.

Birkhoff's theorem can be quickly "proved" by contradiction if the trajectory of a point $x$ is broken into finite pieces in which the mean value is greater than the integral of the function being averaged.
But here, as in Vitali's example, we are dealing with non-measurable objects, since in each orbit we first select a single point and then construct a partition of the orbit. By sacrificing a set of low measure, we implement this idea measurably. The full Rokhlin-Halmos lemma is not required, but it is convenient, and we will use it.

\vspace{2mm}
\bf Theorem 3.2.\it If $\int_X f\, d\mu=0,$ then for an ergodic automorphism $T$ and a sequence of operators $P_N=\sum_{i=0}^{N-1}T^i$, for a.e. $x$, we have
$P_Nf(x) \to 0$. \rm

\vspace{2mm}
Proof by contradiction. Let
$$ \mu (X_s) >0, \ \ X_s=\{x\in X: \limsup P_Nf(x)>s>0\}.$$
By the $T$-invariance of the set $X_s$ and the ergodicity of the automorphism $T$, we have $\mu(X\setminus X_s)=0$, so we can further assume that $X_s=X$.
The set $\{x,Tx,\dots,T^{L(x)-1}x\}$ is called the $s$-heavy $L(x)$-orbit of the point $x$.
We will show that there exists a set of the form $X\setminus \Delta$ such that it can be represented measurably as the union of $L(x)$-orbits on which the mean value of the function $f$ is greater than $s$. Moreover, the measure of the set $\Delta$ can be arbitrarily small.

For all $x$, the function
$$L(x):= \min \{L>0 \, :\, P_Lf(x)>s\} is defined.$$
Let $R$ be a large number such that the measure $\mu(\,x:L(x)>R\,)$ is negligible.
The Rokhlin-Halmos lemma implies that the complement of an arbitrarily small set can be measurably partitioned into $H$-orbits, where $H\gg R$.
In each $H$-orbit, we find the corresponding $s$-heavy $L(x)$-orbits.
To do this, starting from the point $x$, which is the beginning of the $H$-orbit, we move along the $H$-orbit, finding the point $T^ix$ with minimal $i$, which has an $s$-heavy $L(T^ix)$-orbit.
Set $i'=i+L(T^ix)$, move to the point $T^{i'}x$, and search for a new $s$-heavy $L(T^ix)$-orbit for the minimal new $i\geq i'$. The procedure stops when $H-i<L$. Note that the set $\Delta$ of points that are not in the chosen $s$-heavy $L(x)$-orbits can be made arbitrarily small. Since
$\int_{X\setminus\Delta} f\ d\mu\ > \mu(X\setminus\Delta)\, s,$
by the absolute continuity of the integral and the smallness of the measure of the set $\Delta$, we obtain
$\int_{X} f\ d\mu\ > 0. $ A contradiction. Thus, for almost all $x$, we have:
$$\limsup_N P_Nf(x)\leq 0, \ \ \ \limsup_N -P_Nf(x)\leq 0, \ \ \
\lim_N P_Nf(x)=0.$$
The theorem is proved.

This proof cannot be directly extended to actions of the group $\Z^2$.
Let each point of a very large square be covered by a corresponding small square, being its southwest vertex. We can choose a disjoint (!) union of such squares with a total area
significantly smaller than the area of the larger square. 
This is precisely what forces us to modify the proof of the theorem with one-dimensional time.
\section{Rokhlin's Lemma and Wiener's Theorem for $\Z^d$-Actions}

Denote
$$Q_N=\{z=(z_1,\dots,z_d)\, :\, 1\leq z_1,\dots,z_d\leq N\},$$
$$ {Q_N}B=\uu_{z\in Q_N} T^zB.$$

\vspace{2mm}
\bf Lemma 4.1. \it Let $\{T^z\}$ be a free (ergodic) action of the group $\Z^d$.
For every $N$ there exists ${Q_N}B$ of positive measure. \rm

\vspace{2mm}
We leave the proof of the lemma as an exercise. The ergodicity 
is not necessary here.

\vspace{2mm}
\bf Lemma 4.2. \it Let $\{T^z\}$ be a free ergodic action of $\Z^d$, and let $\delta>0$ and towers ${Q_N}B$ and ${Q_H}B$ be of positive measure. There exists $z$ such that
$$\mu({Q_N}B\cap T_g{Q_N}C)> \mu({Q_N}B) \mu({Q_N}C) -\delta.$$
\rm

\vspace{2mm}
\bf Theorem 4.3. \it  Let $\{T^z\}$ be a free ergodic action of $\Z^d$. For every $N$ and $\eps>0$, there exists a tower ${Q_N}B$ such that $\mu({Q_N}B)>1-\eps$.
\rm

\vspace{2mm}
Proof.  Set $a_N=sup \{\mu({Q_N}B)\}$. If $a_N=1$, everything is proved. Let $a_N<1$ and fix a tower $Q_HB'$ such that $v=\mu({Q_H}B')>0$, $H\gg N$, and $H$ is a multiple of $N$. Consider a tower $U=Q_N B$,
whose measure $u$ is extremely close to $a_N$. By Lemma 4.1, there exists a
tower $V=Q_HT^zB'$ such that
$$\mu(U\cap V) < \, uv +\delta.$$
We subtract from $UQ_N B$ sets of the form $Q^N\{x\}$, $x\in B$, which have a nonempty intersection with tower $V$. We obtain a tower $U'$ of height $N$. The height $H$ of tower $V$ is large compared to the height $N$
of the removed sets. Therefore, the decrease in the measure of tower $U$ does not exceed
$$v-v\left(1-\frac N H\right)^d.$$
Therefore, for every
$\delta>0$, we find a tower $U'$ for which $\mu(U')>u-\delta$.
Since $H$ is a multiple of $N$, tower $V$ can also be viewed as a tower of height $N$ with a different base. Combining $V$ with tower $U'$ (they are disjoint) yields a tower $\tilde U$ of height $N$. Note that we have added significantly more to $U'$ than we have subtracted from $U$. This leads to a contradiction.
Indeed, since $a_N<1$, for sufficiently small $\delta$ we obtain
$$a_N\geq \mu(\tilde U) > a_N -\delta +v - (1-a_N +\delta)v -2\delta\, > a_N.$$
Thus, $a_N=1$, the theorem is proved.

For a more general Rokhlin lemma for amenable group actions, see \cite{OW}. Now we are ready to prove almost-everywhere convergence for $\Z^d$-actions.

\vspace{2mm}
\bf Theorem 4.4. \it If $\int_X f\, d\mu=0,$ then for an ergodic free $\Z^d$-action $\{T^z\}$ we have a.e. 
$$P_Nf(x)= N^{-d} \sum_{z\in Q_N}T^zf(x)\to 0, \ \ N\to\infty.$$ \rm

\vspace{2mm}
Proof. Let $ \mu (X_s) >0, \ \ X_s=\{x\in X: \limsup P_Nf(x)>s>0.$
By virtue of the $T^z$-invariance of the set $X_s$ and the ergodicity of our action,
without loss of generality, we assume that $X_s=X$.
The set $T^{Q_{L(x)}}$ is called the $s$-heavy $L(x)$-orbit of the point $x$.
We choose a large number $L$ and $N\gg L$ such that for most $x\in X$,
there is $L(x)\in [L,N]$. We consider a tower $U=T^{Q_H}B$ of measure very close to 1,
where $H\gg N$. Most points $x$ in $T^{Q_H}B$ are covered by squares
$Q_x= T^{Q_{L(x)}}\{x\}$. Consequently, for most square orbits of size $H$,
which make up $U$, most orbital points in them
are covered by squares $Q_x$. There exists a positive constant $c_d$ such that
from any covering of the large cube by small cubes, one can choose
a family of disjoint small cubes (with which they were covered) such that its volume
is no less than the $c_d$-fraction of the volume of the large cube.

If this is not true, then the maximum volume of the optimal family will be very small,
and will become small if we greatly increase the sizes of the cubes in the family. Then there will be many points in the large cube that are far removed from the cubes in the family. But they can be covered by corresponding cubes so that the measure of the optimal covering increases, which contradicts optimality (we leave the details as an exercise).

We define the set $Y_N\subset X$ as the union of all points within the optimal
coverings. Asymptotically, the measures of the sets $Y_N$ are not less than the number $c_d$, and by the definition of these sets, we have
$\limsup_N\int_{Y_N}f d\mu \geq \, c_d.$
The sets $Y_N$ are such that $\mu(Y_N \Delta T^zY_N)\to 0$,
so, by use of   Theorem 2.2. we see that
$$ \lim_N \frac 1 {\mu( Y_N)} \int_{Y_N}f d\mu \to \int_{X}f d\mu=0.$$
We get  a contradiction. The rest is obvious. The theorem is proven.

\section
{Value distributions of ergodic averages}
Recall that the distribution of the values of a function $f:X\to \R$ is understood as the projection $\pi(f)$  onto $\R$ of a measure that is the lift onto the graph $\{(x,f(x)):x\in X\}$  the probability measure $\mu$ on $X$. The closeness of distributions is understood as the closeness of the corresponding functionals on $C_0(\R)$ in the $\ast$-weak sense.

Let an ergodic free action be fixed with  a sequence of functions  $d_j:X\to\R$,
$\|d_j\|_1=1. $
We are looking for a function $f$ and a rapidly growing sequence $N_j$ such that the distribution of the averages $P_{N_j}f(x)/\|P_{N_j}f/|_1$ 
will be  close as we wish  to  $\pi(d_j)$.

For this we consider a tower of  height $jh_j$, then  divide it into $j^d$ subtowers of height $h_j$. We define a function $f_j$ that takes an arbitrary constant value on each subtower. We assume that the positive measure of the complement $E_j$ of the tower is much less than $j^{-j}$. We define the constant value of $f_j$ on $E_j$ so that
over the entire space, $f_j$ has zero mean. We want to have  
$$\pi \left(\frac {P_{n_j}f}   {\|P_{n_j}f\|}\right)\approx    \pi(d_j).$$
Let's act in the spirit of the works  \cite{25}, \cite{PR}. Consider a series $f= \sum_j f_j$ such that the norms $\|f_j\|$ rapidly decrease, and the functions $\sum_{k=1}^{j-1} P_{N_j}f_k$ are extremely small compared to $P_{N_j}f_j\approx f_j$. This is possible because $P_{N}f_k$ for fixed $f_k$, $k<j$, tend to zero and, starting from $N_j$, are as small as we need.
So we choose $N_j$. We make the remainder $\Delta_j = \sum_{m=j+1}^{\infty} P_{N_j}f_m$ very  small compared to
$P_{N_j}f_j$ by choosing very very little functions $f_m$, $m>j$, and therefore the function
$P_{N_j}\Delta_j $ are very small compared to $f_j$.
It turns out that the distributions of the values of the functions
$$ d_j, \ \ \frac {P_{N_j}f }{\|P_{N_j}f\|_1 }, \ \ \frac {P_{N_j}f_j }
 {\|P_{N_j}f_j\|_1 }$$
differ from each other arbitrarily little.
For our  $d_j$ we define $f_j$ ofzero  integral, note that this may have little effect on the distribution of values, since for this we can change the function over a set of a very small measure.

Thus, we have learned to control the value distributions  of the functions $\frac {P_{N_j}f }{\|P_{N_j}f\|_1 }$ for a very sparse sequence $N_j$.
However, it is easy to make it very dense. Since the $f_j$ are constant on almost invariant sets, it is easy to ensure the following:
for $n\in [N_k, R_k]$, all $\frac {P_{n}f }{\|P_{n}f\|_1 }$ differ little from each other, and the ratio $R_k/N_k$ can tend to infinity at any rate.
Note that in this case, the ratio $N_{k+1}/N_k$ should tend to infinity significantly faster, but here we are not constrained in any way in constructing our function
$f$.

\bf Remark. \rm In connection with a question in \cite{P} on optimal estimates for the rate of the  a.e.-convergence we note another property of our averages. This property is an possible  asymptotic independence of the functions
$$\frac {P_{N_k}f }{\|P_{N_k}f\|_1 }, \ \ \frac {P_{N_j}f }{\|P_{N_j}f\|_1 }, \ k\neq j.$$ One can   show  with it   a  lack of reasonable estimates for the rate of the  a.e.-convergence  of  $P_Nf$.  But this requires separate consideration.

\vspace{5mm}
The author thanks Jean-Paul Thouvenot and Benjamin Weiss 
for discussions and comments.

\large


\begin{thebibliography}{9}
\bibitem{K} U. Krengel, Ergodic Theorems, De Gruyter Stud. Math., 6, de Gruyter, Berlin, 1985

\bibitem{T}T. Eisner, B. Farkas, A Journey Through Ergodic Theorems, Birkhauser Advanced Texts Basler Lehrbucher, Birkhauser Verlag, 2025

\bibitem{KSF} I. P. Kornfeld, Ya. G. Sinai, S. V. Fomin, Ergodic Theory, Nauka, Moscow, 1980

\bibitem{BR} I.V. Bychkov, V.V. Ryzhikov,
Singular geometric averages for ergodic multiflows, arXiv:2605.12695 (to appear in
Math. Notes)

\bibitem{FKW}  H. Furstenberg, Y. Katznelson, B. Weiss, Ergodic theory and configurations in sets of positive density
Mathematics of Ramsey theory 5, 184-198

\bibitem{OW} D. S. Ornstein, B. Weiss, Entropy and isomorphism theorems for actions of amenable groups. J. Anal. Math. 48 (1987), 1-141

\bibitem{25}V.V. Ryzhikov, Slow convergence of weighted averages for flows and actions of countable amenable groups, Russian Math. Surveys, 80:5 (2025), 915-918  

\bibitem{PR}  I.V. Podvigin, V.V. Ryzhikov, The full diapason of convergence rates of Birkhoff averages for ergodic flows. arXiv:2601.21404 (to appear in Siberian Math. J.)

\bibitem{P} I.V. Podvigin,   On convergence rates in the Birkhoff Ergodic Theorem, Siberian Math. J., 65:5 (2024), 1170-1186

\end{thebibliography}
\end{document}